\documentclass[12pt]{amsart}

\newtheorem{theorem}{Theorem}
\newtheorem*{corollary}{Corollary}
\newtheorem*{lemma}{Lemma}

\newcommand{\F}{\mathcal{F}}
\newcommand{\M}{\mathcal{M}}
\newcommand{\C}{\mathbb{C}}

\DeclareMathOperator{\re}{Re}
\newcommand{\sumprime}{\sideset{}{'}{\sum}}
\newcommand{\prodprime}{\sideset{}{'}{\prod}}

\begin{document}

\title{Majorant Series}

\author{Harold P. Boas}

\address{Department of Mathematics, Texas A\&M University,
  College Station, TX 77843--3368, USA}

\email{boas@math.tamu.edu}

\thanks{This article is based on a lecture at the third Korean
  several complex variables symposium in December 1998.  The
  author thanks GARC, the Global Analysis Research Center at
  Seoul National University, for sponsoring his participation in
  this international conference.}  

\thanks{The author's research was partially supported by grant
  number DMS-9801539 from the National Science Foundation of the
  United States of America.}

\begin{abstract}
  This article discusses questions in one and several complex
  variables about the size of the sum of the moduli of the terms
  of the series expansion of a bounded holomorphic function.
  Although the article is partly expository, it also includes
  some previously unpublished results with complete proofs.
\end{abstract}


\maketitle

%% modify AMS default to indent subsection in table of contents
\makeatletter
\renewcommand{\l@subsection}{\@tocline{2}{0pt}{3pc}{5pc}{}}
\makeatother

\tableofcontents

\section{Introduction}
What properties of a holomorphic function can be detected from
the \emph{moduli} of its Maclaurin series coefficients?

To make this question precise, fix a holomorphic function~\(f\)
such that \(f(z)=\sum_{k} c_{k} z^{k}\), and consider the class
$\F$ of all power series expansions $\sum_{k} b_{k} z^{k}$ with
the property that $|b_{k}| = |c_{k}|$ for all~$k$. Which
properties of~$f$ are inherited by all functions in the
class~$\F$?

Here are some examples of properties that hold for all members
of~\(\F\) if they hold for one member of~\(\F\).
\begin{itemize}
\item The radius of convergence of the Maclaurin series
  equals~\(1\) (since Hadamard's formula for the radius of
  convergence depends only on the moduli of the coefficients).
\item The function belongs to the Hardy space \(H^2\) of the unit
  disk (since the Hardy space norm equals the square root of the
  sum of the squares of the moduli of the Maclaurin series
  coefficients).
\item The function is univalent in a neighborhood of
    the origin (since this property holds if and only if
    \(c_1\ne0\).)
\end{itemize}

On the other hand, here are some examples of properties that may
hold for some members of~\(\F\) but not for others.
\begin{itemize}
\item The function is holomorphic and univalent in the whole unit
  disk. 
  
  For instance, the linear fractional transformation
  \(z\mapsto\frac{z}{1-z}\) is univalent in the unit disk, but
  changing the sign of the first term of the Maclaurin series
  yields the function \(z\mapsto -2z + \frac{z}{1-z}\), which
  maps the points \(0\) and~\(1/2\) both to~\(0\).

\item The function has zeroes at prescribed locations. 
  
  Indeed, composing a member of~\(\F\) with the reflection
  \(z\mapsto-z\) yields a new member of~\(\F\) with zeroes at
  different locations.

\item The function is holomorphic and bounded in the unit disk.
  
  For instance, for almost every choice of plus and minus signs,
  the series \(\sum_{k=1}^\infty \pm z^k/k\) is continuous when
  \(|z|\le1\) (see, for example, \cite[Chapter~V,
  Theorem~8.34]{Zygmund}); but if all plus signs are taken, the
  series is unbounded in the unit disk.
\end{itemize}

\section{Bohr's theorem}
More generally, it might happen that if $f$~has a certain
property, then every function in the associated class~$\F$ has
some related property.  For example, an old theorem of Harald
Bohr \cite{BohrH1914} implies that if \(f\)~is in the unit ball
of \(H^\infty\) of the unit disk, then every element of the
class~\(\F\) is in the unit ball of \(H^\infty\) of the disk of
radius~\(1/3\).

\begin{theorem}[Bohr, 1914]
  If \(|\sum_{k=0}^\infty c_k z^k|< 1\) when \(|z|< 1\), then
  \(\sum_{k=0}^\infty |c_k z^k|< 1\) when \(|z|< 1/3\). Moreover,
  the radius \(1/3\) is the best possible.
\end{theorem}

This surprising theorem has been largely forgotten. The following
proof, based on a classical inequality of Carath\'eodory, is due
to Edmund Landau \cite{LandauE}.

\begin{lemma}[Carath\'eodory's inequality]
  If \(g\)~is a holomorphic function with positive real part in
  the unit disk, and \(g(z)=\sum_{k=0}^\infty b_k z^k\), then
  \(|b_k|\le 2 \re b_0\) when \(k\ge1\).
\end{lemma}

\begin{proof}[Proof of Bohr's theorem]
  Let \(f(z)\) denote the series \(\sum_{k=0}^\infty c_k z^k\).
  Let \(\varphi\) be an arbitrary real number, and set the
  function~\(g\) in Carath\'eodory's inequality equal to \(1-
  e^{i\varphi}f\) to deduce that \(|c_k|\le 2 \re
  (1-e^{i\varphi}c_0)\) when \(k\ge1\). Since \(\varphi\)~is
  arbitrary, it follows that \(|c_k|\le 2(1-|c_0|)\) when
  \(k\ge1\).
  
  If \(f\)~is a constant function, then the conclusion of the
  theorem is trivial. For nonconstant~\(f\), the preceding
  inequality shows that if \(|z|<1/3\), then
\begin{equation*}
  \sum_{k=0}^\infty |c_k z^k| < |c_0| + 2(1-|c_0|)
  \sum_{k=1}^\infty (1/3)^k = 1.
\end{equation*}

To see that the radius~\(1/3\) in Bohr's theorem is the best
possible, consider the linear fractional transformation~\(f_a\)
defined by \(\displaystyle f_a(z)=\frac{z-a}{1-az}\) when
\(0<a<1\).  Writing \(f_a(z)= \sum_{k=0}^\infty c_k(a) z^k\), one
easily computes that \(\sum_{k=0}^\infty |c_k(a) z^k| = 2a +
f_a(|z|)\).  Simple algebra shows that \(2a + f_a(|z|) >1\) when
\(|z|>1/(1+2a)\). Since \(a\)~can approach~\(1\) from below, it
follows that the radius~\(1/3\) in Bohr's theorem cannot be
increased.
\end{proof}

\begin{proof}[Proof of Carath\'eodory's inequality]
  A common proof of Carath\'eodory's inequality (see
  \cite[page~41]{DurenPL}, for example) uses the Herglotz
  representation for positive harmonic functions. I learned the
  following even simpler proof of the inequality from
  \cite{AizenbergAytunaDjakov}.
  
  By considering \(g(rz)\) and letting \(r\)~increase
  toward~\(1\), we may assume without loss of generality that
  \(g\)~is holomorphic in a neighborhood of the closed unit disk.
  When \(k\ge1\), orthogonality implies that
  \begin{equation*}
    b_k = \frac{1}{2\pi}\int_{0}^{2\pi} e^{-ik\theta}
    g(e^{i\theta}) \,d\theta = \frac{1}{2\pi}\int_{0}^{2\pi}
    e^{-ik\theta}
    \left(
      g(e^{i\theta}) + \overline{g(e^{i\theta})}
    \right)
     \,d\theta.
  \end{equation*}
Consequently,
\begin{equation*}
  |b_k| \le \frac{1}{2\pi}\int_{0}^{2\pi}
  |  2\re g(e^{i\theta}) | \,d\theta = 2 \re b_0,
\end{equation*}
where the last step follows because \(2\re g\) is a positive
function with the mean-value property.
\end{proof}

A natural way to study the class~\(\F\) associated to a
holomorphic function~\(f\) is to single out a canonical
representative of the class. If \(f(z)=\sum_{k=0}^\infty c_k
z^k\), then the \emph{majorant function} \(\M f\) is defined by 
\(\M f(z) = \sum_{k=0}^\infty |c_k| z^k\).

Evidently \( \sup_{|z|<r} |f(z)| \le \sup_{|z|<r}
  |\M f(z)|=\M f(r)\), which justifies the name ``majorant
  function''. Bohr's theorem may be restated as the inequality
\(\sup_{|z|<r} |\M f(z)| \le \sup_{|z|<3r} |f(z)|\) for
every~\(r\) for which the right-hand side makes sense.

\section{Wintner's problem}
In 1956, Aurel Wintner \cite{WintnerA1956} raised the following
question related to Bohr's theorem on majorant series:
\begin{quote}
  For the class of holomorphic functions~\(f\) in the unit disk
  with modulus bounded by~\(1\), find the best upper bound on
  \(\displaystyle \inf_{0<|z|<1} \left|\frac{\M f(z)}{z}\right|\)
  or on \(\displaystyle \inf_{0<r<1} \frac{\M f(r)}{r}\).
\end{quote}
It follows by taking \(r\)~equal to~\(1/3\) that the value~\(3\)
is an upper bound. Wintner claimed---incorrectly---that the
bound~\(3\) cannot be improved.

In a subsequent paper \cite{SchlenstedtG1962} whose title quotes
Wintner's title, G\"unther Schlenstedt pointed out that Wintner
made a blunder in high-school algebra, and that it is easy to see
from the maximum principle that the best bound on \(\displaystyle
\inf_{0<|z|<1} \left|\frac{\M f(z)}{z}\right|\) is
actually~\(1\).  It is amusing to note that the address printed
at the end of Schlenstedt's paper is
``Carl-Friedrich-von-Siemens-Schule, Berlin'', so Schlenstedt was
presumably either a student or a teacher at a German high school.

Perhaps one should not be too critical of Wintner for this error,
for he published 26 papers in 1956. Moreover, even the mistakes
of a good mathematician are interesting.

In my view, the really interesting mistake in Wintner's paper is
one that Schlenstedt did not address. Namely, Wintner claimed to
study \(\displaystyle\inf_{0<|z|<1} \left|\frac{\M
    f(z)}{z}\right|\), but he actually studied \(\displaystyle
\inf_{0<r<1} \frac{\M f(r)}{r}\), apparently thinking that these
quantities admit the same optimal bound. In fact, the best bound
on the second quantity is neither \(3\) nor~\(1\), but~\(2\), as
I shall now demonstrate.

\begin{theorem}
  If \(f\)~is a holomorphic function such that \(|f(z)|<1\) when
  \(|z|<1\), then the majorant function \(\M f\) satisfies the
  inequality
  \begin{equation*}
    \inf_{0<r<1} \frac{\M f(r)}{r} \le2,
  \end{equation*}
  and the bound~\(2\) cannot be replaced by any smaller number.
\end{theorem}

\begin{proof}
  I shall prove somewhat more than is stated in the
  theorem. Namely, I do not need the function~\(f\) to belong to
  the unit ball of \(H^\infty\); all I will use is that \(f\)
  belongs to the unit ball of the Hardy space~\(H^2\). In other
  words, a sufficient hypothesis on~\(f\) is that 
  \begin{equation*}
    \sup_{0<r<1} \frac{1}{2\pi} \int_0^{2\pi} |f(r
    e^{i\theta})|^2 \,d\theta \le1.
  \end{equation*}
  
  By Parseval's theorem, this hypothesis implies that if \(f(z)\)
  has the series expansion \(\sum_{k=0}^\infty c_k z^k\) when
  \(|z|<1\), then \(\sum_{k=0}^\infty |c_k|^2\le1\).
  Consequently, by the Cauchy-Schwarz inequality,
\begin{align*}
  \frac{\M f(r)}{r} &= \frac{|c_0|}{r} + \sum_{k=1}^\infty |c_k|
  r^{k-1} \le \frac{|c_0|}{r} + \biggl( \sum_{k=1}^\infty |c_k|^2
  \biggr)^{1/2} \frac{1}{\sqrt{1-r^2}} \\
&\le \frac{|c_0|}{r} + \frac{\sqrt{1-|c_0|^2}}{\sqrt{1-r^2}}.
\end{align*}
Let \(r\to|c_0|\) to deduce that \(\displaystyle \inf_{0<r<1}
\frac{\M f(r)}{r} \le2\).

For the particular function~\(f\) such that \(\displaystyle
f(z)=\frac{z-\frac{1}{\sqrt{2}}}{1-\frac{z}{\sqrt{2}}}\), a
straightforward computation shows that
\(\displaystyle\inf_{0<r<1}\frac{\M f(r)}{r}=2\), so the
bound~\(2\) cannot be improved.
\end{proof}

\section{A multi-dimensional analogue of Bohr's theorem}
Some researchers in several complex variables feel that power
series are boring, because the really interesting parts of
multi-dimensional complex analysis are the parts that differ from
the one-variable theory, while the elementary theory of power
series superficially appears the same in all dimensions. My goal
here is to exhibit some interesting and accessible problems about
multi-variable power series in which the dependence on the
dimension is the key issue.

I shall use multi-index notation to write an \(n\)-variable power
series as \(\sum_\alpha c_\alpha z^\alpha\), where
\(\alpha\)~denotes an \(n\)-tuple \((\alpha_1, \dots, \alpha_n)\)
of non-negative integers, and \(z^\alpha\)~denotes the product
\(z_1^{\alpha_1}\dotsm z_n^{\alpha_n}\). It is standard (see
\cite[pages 78--80]{RangeRM}, for example) that such a power
series converges in a logarithmically convex complete Reinhardt
domain. To each such domain~\(G\) corresponds a \emph{Bohr
  radius} \(K(G)\): the largest~\(r\) such that whenever
\(|\sum_{\alpha} c_\alpha z^\alpha|\le1\) for \(z\) in~\(G\), it
follows that \(\sum_\alpha |c_\alpha z^\alpha|\le1\) for \(z\) in
the scaled domain~\(rG\).

The terminology ``Bohr radius'' is somewhat whimsical, for
physicists consider the Bohr radius~\(a_0\) of the hydrogen atom
to be a fundamental constant: its value is \(4\pi\epsilon_0
\hbar^2/m_e e^2\), or about \(0.529\)~\AA\@.  The physicists'
Bohr radius is named for Niels Bohr, a founder of the quantum
theory and the 1922 recipient of the Nobel Prize for physics.
Since Niels was the elder brother of Harald Bohr, adopting the
term ``Bohr radius'' for mathematical purposes keeps the honor
within the family.

In contrast to the situation in one dimension, there is no
higher-dimensional bounded domain~\(G\) whose Bohr radius
\(K(G)\) is known \emph{exactly}. However, reasonably good bounds
for the Bohr radius are known for special domains like the ball
and the polydisc.

More generally, let \(B_p^n\) denote \(\{\,z\in\C^n: \sum_{j=1}^n
|z_j|^p < 1\,\}\), which is the unit ball of the complex Banach
space~\(\ell_p^n\) whose norm is defined by \(\|z\|_{\ell_p^n}:=
(\sum_{j=1}^n |z_j|^p)^{1/p}\).  The ball \(B_2^n\) corresponding
to the case that \(p=2\) is the usual Euclidean unit ball
in~\(\C^n\). The ball \(B_\infty^n\) is to be interpreted as the
unit polydisc in~\(\C^n\). The following theorem quantifies the
rate at which the Bohr radius of~\(B_p^n\) decays as the
dimension~\(n\) increases.

\begin{theorem}
\label{thm:multi}
When \(n>1\), the Bohr radius \(K(B_p^n)\) of the \(\ell_p^n\)
unit ball in~\(\C^n\) admits the following bounds.
  \begin{itemize}
  \item If \(1\le p\le 2\), then
    \begin{equation*}
      \frac{1}{3\sqrt[3]{e}} \cdot 
\left(\frac{1}{n}\right)^{1-\frac{1}{p}} \le K(B^n_p)
    < 3\cdot\left(\frac{\log n}{n} \right)^{1-\frac{1}{p}}. 
    \end{equation*}
  \item If \(2\le p\le\infty\), then
    \begin{equation*}
      \frac{1}{3} \cdot \sqrt{\frac{1}{n}} \le K(B^n_p) <
  2\cdot\sqrt{\frac{\log n}{n}}.
    \end{equation*}
  \end{itemize}
\end{theorem}

The theorem does not quite fix the sharp decay rate of the Bohr
radius with the dimension, for the upper bounds contain a
logarithmic factor not present in the lower bounds.  This
logarithmic factor, an artifact of the proof, presumably should
not really be present.

The numerical values of the constants in the theorem can be
improved. For example, the constants in the two parts of the
theorem could be made to agree when \(p=2\).  I have chosen to
write simple constants for clarity, since the main point of the
theorem is the dependence of the Bohr radius on the
dimension~\(n\).

The lower bound in the theorem is due to Lev Aizenberg
\cite{AizenbergL1999} when \(p=1\). The lower bound when \(2\le
p\le\infty\) and the upper bound when \(p=\infty\) are due
jointly to Dmitry Khavinson and myself \cite{BoasKhavinson1997}.
The other parts of the theorem are new.

\section{Proof of the estimates for the Bohr radius}
\subsection{The lower bound}
Lower bounds on the Bohr radius follow from upper bounds on the
Maclaurin series coefficients of a bounded holomorphic function.
In the writings of Bohr~\cite{BohrH1914} and
Landau~\cite{LandauE}, one finds a useful trick employed by
F.~Wiener in the context of coefficient bounds for the
one-dimensional Bohr theorem.

\begin{lemma}[after F.~Wiener]
  Let \(G\) be a complete Reinhardt domain, and let \(\F\) be the
  set of holomorphic functions on~\(G\) with modulus bounded
  by~\(1\). Fix a multi-index~\(\alpha\) other than
  \((0,\dots,0)\), and suppose that the positive real
  number~\(b\) is an upper bound for the modulus of the
  derivative \(f^{(\alpha)}(0)\) for every function~\(f\)
  in~\(\F\).  Then \(|f^{(\alpha)}(0)| \le (1-|f(0)|^2)b\) for
  every function~\(f\) in~\(\F\).
\end{lemma}

\begin{proof}
  Pick an index~\(j\) for which \(\alpha_j\ne0\), and let
  \(\omega_j\)~denote a primitive \(\alpha_j\)th root of unity.
  Consider the function obtained by averaging~\(f\):
  \begin{equation*}
    (z_1, \dots, z_n) \mapsto \frac{1}{\alpha_j} \sum_{k=1}^{\alpha_j}
    f(z_1, \dots, z_{j-1}, \omega_j^k z_j, z_{j+1}, \dots, z_n).
  \end{equation*}
  Iterate this averaging procedure for each non-zero component of
  the multi-index~\(\alpha\). The resulting function~\(h\) is
  still in the class~\(\F\), and its Maclaurin series starts out
  \(f(0) + f^{(\alpha)}(0) z^\alpha/\alpha! + \dotsb\), where
  \(\alpha!\)~denotes the product \(\alpha_1!\dotsm\alpha_n!\).
  Composing~\(h\) with a linear fractional transformation of the
  unit disk that maps \(f(0)\) to~\(0\) gives a new function in
  the class~\(\F\) whose Maclaurin series begins \(
  f^{(\alpha)}(0) z^\alpha/ (1-|f(0)|^2)\alpha! + \dotsb\), and
  the conclusion of the lemma follows.
\end{proof}

To prove the lower bound in Theorem~\ref{thm:multi} when \(1\le
p\le 2\), I will follow the argument used by Aizenberg
\cite{AizenbergL1999} for the case when \(p=1\).  The method
actually yields a bound when \(1\le p\le \infty\), but the bound
is interesting only when \(1\le p\le 2\).

Let \(|\alpha|\) denote the sum \(\alpha_1+\dots+\alpha_n\), and
let \(\alpha^\alpha\) denote the product
\(\alpha_1^{\alpha_1}\dotsm \alpha_n^{\alpha_n}\) (where
\(0^0\)~is interpreted as~\(1\)).  A straightforward calculation
using Lagrange multipliers shows that
\begin{equation}
\label{eq:Lagrange}
\sup\{\,|z^\alpha|: z\in B_p^n\,\} = \left
( \frac{\alpha^\alpha} {|\alpha|^{|\alpha|}}\right)^{1/p}.
\end{equation}
By Cauchy's estimates, it follows that if \(\left| \sum_\alpha
  c_\alpha z^\alpha \right| \le 1\) when \(z\in B_p^n\), then
\(|c_\alpha|\le \left(|\alpha|^{|\alpha|}
  /\alpha^\alpha\right)^{1/p}\le
|\alpha|^{|\alpha|}/\alpha^\alpha\).  Wiener's lemma improves
this estimate by a factor of \((1-|c_0|^2)\). Hence, we have when
\(k>0\) that
\begin{equation*}
\sum_{|\alpha|=k} |c_\alpha z^\alpha| \le (1-|c_0|^2)
\sum_{|\alpha|=k} \frac{k^k}{\alpha^\alpha} |z^\alpha|.
\end{equation*}
Multiplying and dividing by the multinomial coefficient
\(\binom{k}{\alpha}\), which equals \(k!/\alpha!\), observing that 
\(\sum_{|\alpha|=k} \binom{k}{\alpha} |z^\alpha| =
\|z\|_{\ell_1^n}^k\), and applying H\"older's inequality
to bound  \(\|z\|_{\ell_1^n}\) above by \(n^{1-\frac{1}{p}}
\cdot\|z\|_{\ell_p^n}\), we deduce that
\begin{equation*}
\sum_\alpha |c_\alpha z^\alpha| \le |c_0| + (1-|c_0|^2)
\sum_{k=1}^\infty \frac{k^k}{k!} 
\left(n^{1-\frac{1}{p}} \cdot \|z\|_{\ell_p^n}\right)^k.
\end{equation*}
Because the real quadratic function \(t\mapsto t+(1-t^2)/2\)
never exceeds~\(1\), it follows that if \(x\)~is the unique
positive number such that
\begin{equation}
\label{eq:tree}
\sum_{k=1}^\infty \frac{k^k}{k!} x^k =\frac{1}{2},
\end{equation}
then the Bohr radius \(K(B_p^n)\) is at least as big as
\(x/n^{1-\frac{1}{p}}\).

In combinatorics, one encounters the tree function~\(T\) (see,
for example, \cite[p.~395]{KnuthD1997}), which satisfies the
functional equation \(T(x)e^{-T(x)}=x\) and has the series
expansion \(T(x) = \sum_{k=1}^\infty \frac{k^{k-1}}{k!}x^k\).
The equation~\eqref{eq:tree} says that \(xT'(x)=1/2\). Moreover,
the functional equation implies that \(xT'(x)=T(x)/(1-T(x))\), so
\eqref{eq:tree} yields \(T(x)/(1-T(x))=1/2\), or \(T(x)=1/3\).
The functional equation gives the solution
\(x=1/(3\sqrt[3]{e}\,)\).  One can also read off the solution
of~\eqref{eq:tree} from \cite[p.~707]{Prudnikov}, a source
brought to my attention by Lev Aizenberg.  This completes the
proof of the lower bound in Theorem~\ref{thm:multi} when \(1\le
p\le 2\).

A different argument shows that a better lower bound for the Bohr
radius is available when \(2<p\). I reproduce the proof
from~\cite{BoasKhavinson1997}.

If \(|\sum_\alpha c_\alpha z^\alpha|\le 1\) when \(z\in B_p^n\),
then applying the one-dimensional version of Wiener's lemma to
the one-variable function \(\zeta\mapsto \sum_\alpha c_\alpha
\zeta^{|\alpha|} z^\alpha\) shows that \(|\sum_{|\alpha|=k}
c_\alpha z^\alpha| \le 1-|c_0|^2\) when \(k\ge1\). Integrating
the square of the left-hand side of this inequality over a torus
and using the orthogonality of the monomials~\(z^\alpha\) shows
that \((\sum_{|\alpha|=k} |c_\alpha|^2 |w^\alpha|^2)^{1/2} \le
1-|c_0|^2\) for every point~\(w\) in the unit ball
of~\(\ell_p^n\). Now if \(z\)~lies in the \(\ell_p^n\)~ball of
radius \(1/(3\sqrt{n}\,)\), then applying the Cauchy-Schwarz
inequality with \(w\)~equal to \(3\sqrt{n}\,z\) shows that
\begin{align*}
  \sum_{|\alpha|=k}|c_\alpha z^\alpha| &\le \biggl(\sum_{|\alpha|=k}
  |c_\alpha|^2 |w^\alpha|^2 \biggr)^{1/2} \biggl(\sum_{|\alpha|=k}
  \left(\frac{1}{3\sqrt{n}} \right)^{2k}\biggr)^{1/2}\\
  &\le (1-|c_0|^2)/3^k,
\end{align*}
since \(\sum_{|\alpha|=k}1\le n^k\).
Consequently, we have for such points~\(z\) that
\begin{equation*}
  \sum_\alpha |c_\alpha z^\alpha| \le |c_0| +
  (1-|c_0|^2)\sum_{k=1}^\infty 3^{-k} = |c_0| + (1-|c_0|^2)/2 \le 1.
\end{equation*}
This proves the lower bound in Theorem~\ref{thm:multi} when \(2\le
p\le \infty\).

\subsection{The upper bound}
The case of the polydisc, which corresponds to \(p\)~equal
to~\(\infty\), is considered in~\cite{BoasKhavinson1997}. The
idea is to use probabilistic methods to construct a homogeneous
polynomial having relatively small supremum but relatively large
majorant function. The argument in~\cite{BoasKhavinson1997}
applies a result from~\cite{Kahane1985} on random trigonometric
polynomials.

To handle the case of general~\(p\), I will use a related random
technique from~\cite{ManteroTonge1980}. Unfortunately, the
estimate I require is not stated explicitly in that paper: for
their purposes, the authors did not need to keep track of the
constants in the proof, and besides their argument has to be
modified for the case of random tensors that are symmetric.
Consequently, I repeat here the argument
from~\cite{ManteroTonge1980} with appropriate modifications and
amplifications.
  
Although the goal is to show the existence of a homogeneous
polynomial of degree~\(d\) in \(n\)~variables that satisfies
suitable estimates, technical considerations suggest proving a
more general result for symmetric multi-linear functions mapping
\( (\C^n)^d\) into~\(\C\).

\begin{theorem}
\label{thm:form}
If \(1\le p\le
\infty\), and if \(n\) and~\(d\) are integers larger than~\(1\),
then there exists a symmetric multi-linear function \(F:
(\C^n)^d\to\C\) of the form
  \begin{equation*}
    F(Z_1, \dots, Z_d)= \sum_{J_1=1}^n \dots\sum_{J_d=1}^n \pm
    Z_{1J_1}\dotsm Z_{dJ_d}
  \end{equation*}
  such that the supremum of \(|F(Z_1, \dots, Z_d)|\) when every
  \(n\)-vector~\(Z_k\) lies in the unit ball of \(\ell_p^n\) is
  at most
\begin{equation}
\label{eq:randombound}
\sqrt{32d\log(6d)}\, \times
  \begin{cases}
    n^{\frac{1}{2}} (d!)^{1-\frac{1}{p}}, 
        & \text{if } 1\le p\le 2;\\
    n^{\frac{1}{2}+(\frac{1}{2}-\frac{1}{p})d} (d!)^{\frac{1}{2}},
          & \text{if } 2\le p\le\infty.
  \end{cases}
\end{equation}
\end{theorem}

In the theorem, the plus and minus signs are chosen independently
for each set of indices \(J_1\), \dots, \(J_d\), with the proviso
that the same choice is made for every permutation of a given set
of indices. In other words, \(F\)~is the sum of the
\(n^d\)~monomials of degree~\(d\) that can be formed from the
components of the vectors \(Z_1\), \dots, \(Z_d\), with plus and
minus signs chosen to enforce symmetry under permutations of the
vectors.

When all the vectors~\(Z_k\) are equal, the theorem provides a
special homogeneous polynomial.  Observe that the number of
\(d\)-tuples \(J_1\), \dots, \(J_d\) that are permutations of a
given list of indices is typically less than~\(d!\) because some
indices in the list will be repeated.  Indeed, if
\(\alpha_k\)~denotes the number of times the integer~\(k\)
appears in the list \(J_1\), \dots, \(J_d\), then the number of
\(d\)-tuples that are permutations of the list is the multinomial
coefficient \(\binom{d}{\alpha}\). Consequently, we have the
following corollary of the theorem.

\begin{corollary}
  If \(1\le p\le \infty\), and if \(n\) and~\(d\) are integers
  larger than~\(1\), then there exists a homogeneous polynomial
  of degree~\(d\) in the variable~\(z\) in~\(\C^n\) of the form
  \begin{equation}
    \sum_{|\alpha|=d} \pm\binom{d}{\alpha} z^\alpha
  \end{equation}
  such that the supremum of the modulus of the polynomial when
  \(z\)~lies in the unit ball of~\(\ell_p^n\) is no greater than
  the above bound~\eqref{eq:randombound}.
\end{corollary}

To prove the upper bound in Theorem~\ref{thm:multi}, consider the
homogeneous polynomial of the corollary when \(z\)~is the vector
\((1, \dots, 1)\) scaled by the factor \(K(B_p^n)/n^{1/p}\).  To
write formulas that work for both ranges of~\(p\) simultaneously,
I will write \(m(p)\) for \(\min(p,2)\) and \(M(p)\) for
\(\max(p,2)\).  The definition of the Bohr radius implies that
\begin{equation}
  \sum_{|\alpha|=d} \binom{d}{\alpha} (K(B_p^n)/n^{1/p})^d
\le   
n^{\frac{1}{2}+(\frac{1}{2}-\frac{1}{M(p)})d}
(d!)^{1-\frac{1}{m(p)}} 
        \sqrt{32d\log(6d)}.
\end{equation}
Since \(\sum_{|\alpha|=d} \binom{d}{\alpha} =n^d\), it follows
that
\begin{equation}
\label{eq:stir}
  K(B_p^n)\le 
\left(  \frac{(d!)^{1/d}}{n}\right)^{1-\frac{1}{m(p)}}
 \left(
    32nd\log(6d)\right)^{\frac{1}{2d}}.
\end{equation}

Stirling's formula implies that if \(d\approx \log n\), then the
right-hand side of~\eqref{eq:stir} yields the upper bound in
Theorem~\ref{thm:multi} for sufficiently large~\(n\).  Numerical
calculations show that if \(d=2+\lfloor \log n\rfloor\), then
\(n\)~greater than~\(148\) is sufficiently large.  On the other hand,
the upper bound in Theorem~\ref{thm:multi} holds automatically for
smaller values of the dimension~\(n\), because the multi-dimensional
Bohr radius does not exceed the one-dimensional Bohr radius
of~\(1/3\), while \(2\sqrt{\log n}/\sqrt{n} >1/3\) when \(1<n<189\).
Thus the upper bound in Theorem~\ref{thm:multi} is a consequence of
the corollary of Theorem~\ref{thm:form}.

\begin{proof}[Proof of Theorem~\ref{thm:form}]
  The proof consists of a probabilistic estimate and a covering
  argument. Three parameters are fixed throughout the proof: the
  dimension~\(n\), the degree~\(d\), and the exponent~\(p\).

I will use the following notation: \(Z\)~denotes a \(d\)-tuple
\(Z_1\), \dots, \(Z_d\) of \(n\)-vectors, \(J\)~denotes a
\(d\)-tuple \(J_1\), \dots, \(J_d\) of integers between \(1\)
and~\(n\), sums and products run over all such \(d\)-tuples of
integers, a prime on a sum or product means that it is restricted
to \(d\)-tuples of integers that are arranged in non-decreasing
order, and \(J\sim K\) indicates that the \(d\)-tuples \(J\)
and~\(K\) are permutations of each other. The symbol~\(Z_J\) is
shorthand for the monomial \(Z_{1J_1}\dotsm Z_{dJ_d}\). Thus, the
function in the statement of Theorem~\ref{thm:form} can be
written in the form
\begin{equation*}
  F(Z)= \sumprime_K \biggl(\pm \sum_{J\sim K} Z_J \biggr),
\end{equation*}
and in this expression, all of the plus and minus signs are
independent of each other.

To begin the probabilistic argument, fix a point~\(Z\)
in~\((\C^n)^d\) such that every~\(Z_k\) lies in the \(\ell_p^n\)
unit ball~\(B_p^n\).  For each \(d\)-tuple~\(K\) in
non-decreasing order, choose a different Rademacher
function~\(r_K\).  (Recall that the Rademacher functions are
independent functions each taking the values \(+1\) and \(-1\)
with probability~\(1/2\).)  Consider the random sum
\begin{equation}
F(t,Z):= \sumprime_K \biggl( r_K(t) \sum_{J\sim K} Z_J \biggr),
\end{equation}
where \(t\) lies in the interval \([0,1]\). The immediate goal is
to make an upper estimate on the probability that this sum has
large modulus.

Let \(\lambda\)~be an arbitrary positive real number;
in~\eqref{eq:parameters} below, I will specify a value
for~\(\lambda\) in terms of \(n\), \(d\), and~\(p\).  Invoking
the independence of the Rademacher functions, we can compute the
expectation (that is, the integral with respect to~\(t\)) of the
exponential of the real part of \(\lambda F(t,Z)\) by computing
the product over non-decreasing \(d\)-tuples~\(K\) of the
expectations of \( \exp( \lambda r_K(t) \sum_{J\sim K} \re
Z_J)\): namely,
\begin{equation}
\prodprime_K \cosh \biggl(\lambda \sum_{J\sim K} \re Z_J\biggr).
\end{equation}

In view of the inequality \(\cosh x\le \exp(x^2/2)\), this
expectation is bounded above by
\begin{equation}
\label{eq:expect}
\exp\biggl( \frac{1}{2}\lambda^2 
\sumprime_K \biggl(\sum_{J\sim K} \re Z_J\biggr)^2 \biggr).
\end{equation}
I will again write \(m(p)\) for \(\min(p,2)\) and \(M(p)\) for
\(\max(p,2)\) in order to analyze the cases \(1\le p\le2\) and
\(2\le p\le\infty\) simultaneously.  H\"older's inequality
implies that
\begin{equation*}
  \biggl(\sum_{J\sim K} \re Z_J\biggr)^2 \le (d!)^{2(1-1/m(p))}
  \biggl(\sum_{J\sim K} |Z_J|^{m(p)} \biggr)^{2/m(p)}.
\end{equation*}
The exponent \(2/m(p)\) is equal to~\(1\) when \(2\le
p\le\infty\), and replacing it by~\(1\) when \(1\le p<2\) can
only increase the right-hand side when each \(n\)-vector~\(Z_k\)
lies in~\(B_p^n\). Consequently, \eqref{eq:expect} is bounded
above by
\begin{equation}
\label{eq:nexpect}
\exp \biggl( \frac{1}{2}\lambda^2 (d!)^{2(1-1/m(p))} 
\|Z_1\|^{m(p)}_{\ell_{m(p)}^n} \dotsm
\|Z_d\|^{m(p)}_{\ell_{m(p)}^n}\biggr).
\end{equation}
Since \(Z_k\in B_p^n\), H\"older's inequality implies that
\(\|Z_k\|^{m(p)}_{\ell_{m(p)}^n}\) is bounded above by
\(n^{1-\frac{2}{M(p)}}\).  Therefore \eqref{eq:nexpect} is
bounded above by
\begin{equation}
\label{eq:fexpect}
\exp\biggl(  \frac{1}{2}\lambda^2 (d!)^{2(1-1/m(p))} 
n^{(1-2/M(p))d} \biggr).
\end{equation} 

Let \(R\)~be an arbitrary positive real number;
in~\eqref{eq:parameters} below, I will specify a value for~\(R\)
in terms of \(n\), \(d\), and~\(p\).  By Chebyshev's inequality,
the upper bound~\eqref{eq:fexpect} on the expectation of
\(\exp(\lambda \re F(t,Z))\) implies that the measure of the set
of points~\(t\) for which \(\re F(t,Z)\) exceeds~\(R\) is at most
\begin{equation}
  \label{eq:cheby}
  \exp\left(-R\lambda + \frac{1}{2}\lambda^2
(d!)^{2(1-1/m(p))} n^{(1-2/M(p))d}\right).
\end{equation}
Symmetric reasoning gives the same estimate for the probability
that \(\re F(t,Z)\) is less than~\(-R\), so the probability that
\(|\re F(t,Z)|\) exceeds~\(R\) is at most \(2\)~times
\eqref{eq:cheby}.  The same argument applies to the imaginary
part of \(F(t,Z)\).  Consequently, the probability that
\(|F(t,Z)|\) exceeds \(\sqrt{2} \,R\) is at most \(4\)~times
\eqref{eq:cheby}.

This probabilistic estimate holds for an arbitrary but
fixed~\(Z\).  The second part of the proof is a covering argument
to produce an estimate that is uniform in~\(Z\).  The following
lemma is well known (see \cite[p.~7]{Milman}, for example).  The
exponent differs from the statement in~\cite{ManteroTonge1980}
because I am considering the complex \(\ell_p^n\)~space.

\begin{lemma}
  If \(\epsilon\)~is a positive real number, then the unit ball
  of the complex space~\(\ell_p^n\) can be covered by a
  collection of open \(\ell_p^n\) balls of radius~\(\epsilon\),
  the number of balls in the collection not exceeding \( (1 +
  2\epsilon^{-1})^{2n}\) and the centers of the balls lying in
  the closed unit ball of~\(\ell_p^n\).
\end{lemma}

\begin{proof}
  Place arbitrarily in the open \(\ell_p^n\) ball of radius
  \((1+\frac{\epsilon}{2})\) a collection of disjoint open
  \(\ell_p^n\) balls of radius~\(\epsilon/2\). Since the
  Euclidean volume of a ball in \(\ell_p^n\) scales like the
  \((2n)\)th power of its radius, an obvious volume comparison
  shows that the number of disjoint balls cannot exceed the
  \((2n)\)th power of the ratio of radii
  \((1+\frac{\epsilon}{2})/(\epsilon/2)\): namely, \((1 +
  2\epsilon^{-1})^{2n}\).
  
  If the collection of disjoint balls is made maximal, then every
  point of the closed \(\ell_p^n\) unit ball must lie within
  \(\epsilon/2\) of some point of one of the balls in the
  collection. Consequently, the balls with the same centers but
  with radius~\(\epsilon\) must cover the closed unit ball.
\end{proof}

Next we need a simple Lipschitz estimate for \(F(t,Z)\). Suppose
that \(Z\) and~\(W\) are points of \((\C^n)^d\) such that all of
the component \(n\)-vectors \(Z_1\), \dots, \(Z_d\) and \(W_1\),
\dots, \(W_d\) lie in~\(B_p^n\), and \(\|Z_k-W_k\|_{\ell_p^n} \le
\epsilon\) for every~\(k\). The multi-linearity of~\(F\) implies
that
\begin{multline*}
F(t, Z_1, \dots, Z_d)= F(t, Z_1-W_1, Z_2, \dots, Z_d)\\ +
F(t, W_1, Z_2-W_2, Z_3, \dots, Z_d) \\
+ F(t, W_1, W_2, Z_3-W_3, Z_4, \dots, Z_d) + \dotsb  \\
+ F(t, W_1, W_2, \dots, W_{d-1}, Z_d-W_d) + F(t, W_1, \dots, W_d).
\end{multline*}
Consequently, \(|F(t,Z)-F(t,W)|\) is at most \(d\) times
\(\epsilon\) times the supremum of \(|F(t,\cdot)|\) over
\((B_p^n)^d\). Taking \(\epsilon\) to be \(1/(2d)\), we see by
the lemma that there is a collection of at most \( (1+4d)^{2nd}\)
points of \((B_p^n)^d\) such that the supremum of
\(|F(t,\cdot)|\) over \((B_p^n)^d\) is no more than twice the
supremum of \(|F(t,\cdot)|\) over this finite collection of
points. 

Applying the preceding probabilistic estimate to each member of
the finite collection of points shows that the probability that
the supremum of \(|F(t,\cdot)|\) over \((B_p^n)^d\) exceeds
\(2\sqrt{2}\,R\) is at most
\begin{equation*}
4(1+4d)^{2nd}\exp\biggl(-R\lambda +
\frac{1}{2}\lambda^2 (d!)^{2(1-1/m(p))} n^{(1-2/M(p))d}\biggr).
\end{equation*}

Now take the following values for the parameters \(R\)
and~\(\lambda\):
\begin{equation}
\label{eq:parameters}
\begin{aligned}
R&:=\left( 2(d!)^{2(1-1/m(p))} n^{(1-2/M(p))d}
\log( 8(1+4d)^{2nd})\right)^{1/2},\\
\lambda&:=\frac{R}{(d!)^{2(1-1/m(p))} n^{(1-2/M(p))d}}.
\end{aligned}
\end{equation}
With these choices, we find that the probability that the
supremum of \(|F(t,\cdot)|\) over \((B_p^n)^d\) exceeds
\(2\sqrt{2}\,R\) is at most~\(1/2\), so we are sure that there
exists a particular value of~\(t\) such that the supremum of
\(|F(t,\cdot)|\) over \((B_p^n)^d\) is no more than
\begin{equation}
\label{eq:finalprob}
\left( 16(d!)^{2(1-1/m(p))} n^{(1-2/M(p))d}
\log( 8(1+4d)^{2nd})\right)^{1/2}.
\end{equation}
The values of the Rademacher functions at this particular value
of~\(t\) produce the pattern of plus and minus signs indicated in
the statement of Theorem~\ref{thm:form}.  Moreover,
\(8(1+4d)^{2nd} < (6d)^{2nd}\) when \(n\) and~\(d\) are both at
least~\(2\), so the upper bound~\eqref{eq:finalprob} is even
smaller than the bound stated in Theorem~\ref{thm:form}.
\end{proof}

\section{Another analogue of Bohr's theorem}
In \cite{AizenbergL1999}, Aizenberg introduced a second Bohr
radius \(B(G)\): namely, the largest~\(r\) such that whenever
\(|\sum_{\alpha} c_\alpha z^\alpha| \le1\) for \(z\in G\), it
follows that \( \sum_\alpha \sup_{z\in rG} |c_\alpha z^\alpha|
\le1\).  Clearly \(B(G)\le K(G)\), because the supremum of a sum
is no bigger than the sum of the suprema.  Equality holds for
polydiscs, because points on the torus have the property that
they maximize \(|z^\alpha|\) for all indices~\(\alpha\)
simultaneously.

The following statement is an analogue of Theorem~\ref{thm:multi}
for the second Bohr radius.

\begin{theorem}
\label{thm:second}
When \(n>1\), the second Bohr radius \(B(B_p^n)\) of the
\(\ell_p^n\) unit ball in~\(\C^n\) admits the following bounds.
  \begin{itemize}
  \item If \(1\le p\le 2\), then
  \begin{equation*}
    \frac{1}{3} \cdot \frac{1}{n} <
1-\left(\frac{2}{3}\right)^{\frac{1}{n}}
\le B(B^n_p)
    < 4\cdot\frac{\log n}{n}.
  \end{equation*}
  \item If \(2\le p\le\infty\), then
\begin{equation*}
  \frac{1}{3} \cdot 
\left(\frac{1}{n}\right)^{\frac{1}{2}+\frac{1}{p}} 
\le B(B^n_p) <
  4\cdot\left(\frac{\log n}{n}\right)^{\frac{1}{2}+\frac{1}{p}}.
\end{equation*}
  \end{itemize}
\end{theorem}

\begin{proof}
  The lower bound is due to Aizenberg \cite{AizenbergL1999} when
  \(1\le p\le2\).  Namely, Cauchy's estimates imply that if
  \(|\sum_\alpha c_\alpha z^\alpha|\le1\) in~\(B_p^n\), then
  \(|c_\alpha|\le 1/\sup\{\,|z^\alpha|: z\in B_p^n\,\}\), and
  Wiener's lemma provides an extra factor of \(1-|c_0|^2\).
  Therefore
\begin{equation*}
  \sum_\alpha \sup\{\,|c_\alpha z^\alpha|: 
\|z\|_{\ell_p^n}\le r\,\}  \le |c_0| +
  (1-|c_0|^2) \sum_{\alpha\ne0} r^{|\alpha|},
\end{equation*}
and this quantity will be bounded above by~\(1\) if
\(\sum_{\alpha\ne0} r^{|\alpha|} \le 1/2\). Since
\(\sum_{\alpha\ne0} r^{|\alpha|} = -1 + 1/(1-r)^n\), we deduce
that the second Bohr radius is no smaller than \(1-(2/3)^{1/n}\).
It is easy to see that this quantity exceeds \(1/(3n)\) when
\(n>1\), because the function \(n\mapsto n\cdot( 1-
(2/3)^{1/n})\) is increasing for positive~\(n\) and takes the
value \(1/3\) when \(n=1\).

When \(p=\infty\), the lower bound is contained in
Theorem~\ref{thm:multi}, since the two Bohr radii are the same
for the polydisc.  When \(2<p<\infty\), the lower bound follows
by observing that an \(\ell_\infty^n\) ball of radius
\(1/n^{1/p}\) nests inside~\(B_p^n\), and so \(B(B_p^n)\ge
B(B_\infty^n)/n^{1/p}\).

To prove the upper bound in Theorem~\ref{thm:second}, apply
the definition of the second Bohr radius to the homogeneous
polynomial of the corollary of Theorem~\ref{thm:form}. 
{From}~\eqref{eq:Lagrange}, it follows that
\begin{equation*}
  \sum_{|\alpha|=d} \binom{d}{\alpha}
  \left(\frac{\alpha^\alpha}{d^d}\right)^{\frac{1}{p}}B(B_p^n)^d 
\le
n^{\frac{1}{2} + (\frac{1}{2} - \frac{1}{M(p)})d} 
(d!)^{1-\frac{1}{m(p)}}
 \sqrt{32d\log(6d)}.
\end{equation*}
Since \(\alpha^\alpha\ge1\), and \(\sum_{|\alpha|=d}
\binom{d}{\alpha} =n^d\), we deduce that
\begin{equation*}
  B(B_p^n)\le 
\frac{1}{n^{\frac{1}{2}+\frac{1}{M(p)}}} d^{\frac{1}{p}} 
(d!)^{\frac{1}{d}(1-\frac{1}{m(p)})} 
   (32 nd\log(6d))^{\frac{1}{2d}}.
\end{equation*}
When \(d= \lfloor\log n\rfloor\), the right-hand side
asymptotically gives the desired power of \((\log n)/n\), and
numerical calculations show that the right-hand side is smaller
than the upper bound stated in Theorem~\ref{thm:second} when \(n>
28\).  On the other hand, the upper bound in the theorem holds
automatically for smaller values of~\(n\), because \((4\log
n)/n>1/3\) when \(1<n<46\).
\end{proof}

\section{Directions for future research}
My goal in this article has been not only to demonstrate some
interesting results and techniques in the theory of
multi-variable power series, but also to suggest some avenues for
further investigation. Here are some concrete problems.
\begin{enumerate}
\item Construct examples to show that the logarithmic term is not
  needed in the upper bound in Theorems \ref{thm:multi}
  and~\ref{thm:second}.
\item Find the \emph{exact} value of the Bohr radius for the
  polydisc in dimension~\(2\). 
\item Estimate the Bohr radius for domains of the form
    \begin{equation*}
      |z_1|^{p_1}+ \dots + |z_n|^{p_n}<1.
    \end{equation*}
\item Generalize Wintner's problem to higher dimensions.
\end{enumerate}

\end{document}